\documentclass[10pt,leqno]{article}
 \usepackage{amsmath}
\usepackage{amssymb}
\usepackage{mathrsfs}

 \def\proclaim{\medskip\noindent}\def\endproclaim{\medskip}

 \def\beqlb{\begin{eqnarray}}\def\eeqlb{\end{eqnarray}}
 \def\beqnn{\begin{eqnarray*}}\def\eeqnn{\end{eqnarray*}}

 \def\<{\langle}\def\>{\rangle}
 \def\gfc{\genfrac{}{}{0pt}{}}

 \def\n{~}

 \def\E{\mathbf{E}}\def\P{\mathbf{P}}\def\R{\mathbb{R}}
 \def\DD{\mathscr{D}}\def\FF{\mathscr{F}}\def\LL{\mathscr{L}}

\begin{document}

\centerline{\bf SUPERPROCESSES OVER A STOCHASTIC FLOW }

\centerline{\bf WITH SPATIALLY DEPENDENT BRANCHING$^*$}

\footnote{$^*$ {This work was supported by the National Natural
Science Foundation of China grant no.10721091, 11471254, and the
Multi-Year Research Grant of the University of Macau No.
MYRG2014-00034-FST.}}

\centerline{{ C. DONG}$^\dagger$}\footnote{$^\dagger$ School of Mathematics and Statistics, Xidian University, Xi'an 710126, Shaanxi Province, People's Republic of China. School of
Mathematical Sciences, Beijing Normal University, Beijing 100875,
People's Republic of China (dczdean@mail.bnu.edu.cn).}

\medskip

{\bf Abstract.} This paper considers a generalized model of [G.
Skoulakis and R. J. Adler, {\it Ann. Appl. Probab.}, 11 (2001),
pp.488-543]. We show the existence of superprocesses in a random
medium (flow) with location dependent branching. Technically, we
make use of a duality relation to establish the uniqueness of the
martingale problem and to obtain the moment formulas which
generalize those of [G. Skoulakis and R. J. Adler, {\it Ann. Appl.
Probab.}, 11 (2001), pp.488-543].

\bigskip


\bigskip

 {\bf 1. Motivation and introduction.}
Superprocesses over a stochastic flow are introduced in \cite{SA01},
where the motions of all particles are subject to the common noise
(flow) and the branching of particles is independent of their
motions. The authors of \cite{SA01} used branching particle systems
approximation to prove the existence of superprocesses over a
stochastic flow (flow superprocesses for short). Moreover, they made
detailed use of the approximating processes to establish the moments
of the flow superprocesses. It was mentioned in \cite{SA01} that the
moment formulas could also be obtained by a dual method, which
``would undoubtedly be more elegant" but was not adopted since it
``does not seem to be straightforward'' (\cite[p.497]{SA01}). In
this paper we will consider a generalized model of \cite{SA01}, in
which the branching of particles is location dependent. A similar
model on Polish space was studied in the first part of \cite{P02}.
When considering superprocesses, a martingale problem is usually
inevitable and the duality approach usually plays a key role in
deriving the uniqueness for the martingale problem. To establish
uniqueness, the approach used in \cite{SA01} was to justify the
duality conditions of \cite{DK82} instead of constructing a dual
process. We shall use the latter as in \cite{DLW01} to show
uniqueness and derive moment formulas for the flow superprocess as
well.

In the rest of this section, we give a concise description of our
model; the reader is referred to \cite{SA01} for a more specific one.
 The main results are given as well as proved in the next section.
In the final section, a further extension is provided. Let $N$ be a
positive integer, which varies whenever necessary. Let $E=\R^d$ with
$\Delta$ its infinity and write $\bar{E}=E\cup\{\Delta\}$, the
one-point compactification of $E$. $C_b(E)$ is the space of bounded
continuous real-valued functions. $C_l(E)$ denotes the subspace of
$C_b(E)$ such that its members have limits at infinity. $C^2_l(E)$
stands for the subspace of $C_l(E)$ such that its members have two
continuous derivatives which have limits at infinity. $C^2_b(E)$
consists of the elements in $C_b(E)$ possessing bounded first and
second partial derivatives. $M_F(E)$ is the space of finite Borel
measures on $E$ endowed with weak convergence topology.
$D_{M_F(E)}[0,\infty)$ is the well-known Skorokhod space and the
meaning of $C_{M_F(E)}[0,\infty)$ is obvious. Let $\Rightarrow$ and
$\rightrightarrows$ denote weak convergence and uniform convergence,
respectively. The superscript $+$ attached to a set will mean its
non-negative subset. Write $\mu(h)\equiv \<h,\mu\>$ as the integral
of $h$ with respect to the measure $\mu$. Throughout the paper, let
$\P$ always denote the probability measure for the probability space
involved and $\E$ the corresponding expectation.

 Let
 $I=\{\alpha=(\alpha_0,\alpha_1,\ldots,\alpha_k):k\geq 0,
\alpha_i\in\{1,2,\ldots\}, 0\leq i\leq k\}$ be the family of
multi-indices, setting
$|\alpha|=|(\alpha_0,\alpha_1,\ldots,\alpha_k)|=k$,
$\alpha-1=(\alpha_0,\ldots,\alpha_{|\alpha|-1})$ and
$\alpha|_i=(\alpha_0,\ldots,\alpha_i)$. Let $n=1,2,\ldots$. Suppose
at time zero that $K_n$ (deterministic) particles, located
separately at $x^n_1,\ldots,x^n_{K_n}\in E$, are given. For $t\geq
0$, write $\alpha\sim_nt$ if and only if $|\alpha|/n\leq
t<(1+|\alpha|)/n$ and $\alpha_0\leq K_n$. Each particle in our model
is labeled by a multi-index in $I$. A particle with label $\alpha$
is understood to be born at time $|\alpha|/n$ and to die at
$(1+|\alpha|)/n$ with $N^{\alpha,n}$ offspring reproduced. For
$\alpha\sim_nt$ (between branching), the motion $Y^{\alpha,n}_t$ of
particle $\alpha$ is determined by
  \beqlb\label{1.1}
dY^{\alpha,n}_t=b(Y^{\alpha,n}_t)dt+e(Y^{\alpha,n}_t)dB^{\alpha,n}_t+c(Y^{\alpha,n}_t)dW^n_t,
\quad Y^{\alpha,n}_0=x^n_{\alpha_0},
 \eeqlb
where $b: \R^d\rightarrow \R^d$, $c: \R^d\rightarrow \R^{d\times
m}$, and $e: \R^d\rightarrow \R^{d\times d}$; $W^n$ is an
$\R^m$-valued Brownian motion, random environment (flow),
independent of which is the family $\{B^{\alpha,n}:\alpha_0\leq
K_n\}$ of $E$-valued Brownian motions stopped at time
$t=(|\alpha|+1)/n$. For each $k$, members in
$\{B^{\alpha,n}:\alpha_0\leq K_n,|\alpha|=k\}$ are conditionally
independent given $\sigma\{B^{\alpha,n}:\alpha_0\leq
K_n,|\alpha|<k\}$, and $B^{\alpha,n}_t=B^{\alpha-1,n}_t$ for
$t\leq|\alpha|/n$. Let $k_n=k/n$ and $a_n=1/n$. Define for
$t\in[k_n,k_n+a_n)$ and $k=0,1,\ldots$
 \beqnn
\FF^n_t=\sigma(B^{\alpha,n},N^{\alpha,n}:|\alpha|<k)
\bigvee\bigcap_{r>t}\sigma(W^n_s,B^{\alpha,n}_s:s\leq
r, |\alpha|=k)
 \eeqnn
and
 \beqnn
\bar{\FF}^n_{k_n}=\FF^n_{k_n}\bigvee\sigma(W^n_s,B^{\alpha,n}_s:s\leq
k_n+a_n,|\alpha|=k).
 \eeqnn
Assume that $\{N^{\alpha,n}:|\alpha|=k\}$ are conditionally
independent given $\bar{\FF}^n_{k_n}$, and
 \beqlb\label{1.2}
\left\{
  \begin{array}{ll}
    \E\left(N^{\alpha,n}|\bar{\FF}^n_{k_n}\right)=
        1+\gamma_n(Y^{\alpha,n}_{k_n+a_n})/n=:\beta_n(Y^{\alpha,n}_{k_n+a_n})& \hbox{} \\
    \mbox{Var}\left(N^{\alpha,n}|\bar{\FF}^n_{k_n}\right)=\sigma_n(Y^{\alpha,n}_{k_n+a_n})^2, & \hbox{}
  \end{array}
\right.
 \eeqlb
where $\gamma_n\in C_l(E)$ and $\sigma_n\in C_l(E)^+$.  Now define
  \beqnn
X^n_t(B)=\frac{\mbox{number of particles in $B$ at time $t$}}{n},
 \eeqnn
where $B$ is a Borel subset of $E$. Intuitively, $X^n_t$
characterizes the mass distribution of the particle system at time
$t$.

It is worth pointing out that compared to \cite[p.493]{SA01}, the
different parts in our model are on the one hand the equation
$(\ref{1.1})$, where $e$ is extended to be non-diagonal. On the
other hand the significant difference lies in the branching
mechanism, which is location dependent as indicated in
$(\ref{1.2})$.

Suppose that there exist $p>2$ and $C>0$ such that
  \beqlb\label{1.3}
\E[(N^{\alpha,n})^p]\leq C \mbox{\ for all $\alpha$ and $n$,\ }
  \gamma_n\rightrightarrows \gamma\in
C_l(E)\mbox{\ and\ }\sigma_n\rightrightarrows\sigma\in C_l(E)^+
\mbox{\ as\ } n\to\infty.
 \eeqlb
$\gamma$ is called the {\it drift function} and $\sigma^2$ the {\it
branching variance}. Let $C_\gamma, C_\sigma$ be the constants such
that $|\gamma_n|,|\beta_n|\leq C_\gamma$ and $\sigma_n\leq C_\sigma$
for all $n$.

\medskip
\noindent{\it Remark 1.1}. For each $\gamma\in C_l(E)$ and
$\sigma\in C_l(E)^+$, there exist random variables $\xi_n$ such that
(\ref{1.2}) and (\ref{1.3}) hold (see \cite[p.143]{P02}), and $p$
can be very close to 2.
\medskip

{\bf 2. Continuous spatially dependent branching.} Based on a dual
method,  we shall discuss in this section the existence and moment
properties of a flow superprocess with aforementioned parameters
$\gamma$ and $\sigma$.
\medskip

{\bf Hypotheses~(LU)}
 \begin{itemize}
 \item[(L)] $ |b(x)-b(y)|+\|c(x)-c(y)\|+\|e(x)-e(y)\|\leq
K|x-y|, \quad x,y\in E$.

\item[(U)] $b_i, c_{il},e_{ik}\in C_l^2(E)$, $i,k=1,\ldots,d,
l=1,\ldots,m$, and for any $N\geq 1$ there exists $\lambda_N>0$ such
that
 \beqnn
\sum_{p,q=1}^N\sum_{i,j=1}^d\xi_i^pd_{ij}(x_p,x_q)\xi_j^q
\geq\lambda_N\sum_{p=1}^N\sum_{i=1}^d(\xi_i^p)^2
  \eeqnn
for $x_1,\ldots,x_n\in E$ and
$(\xi_1^1,\ldots,\xi_d^1;\ldots;\xi_1^N,\ldots,\xi_d^N)\in E^N$,
where $d_{ij}(x,y)=\sum_{k=1}^de_{ik}(x)e_{jk}(y)+a_{ij}^{(m)}(x,y)$
with $a_{ij}^{(m)}(x,y)=\sum_{l=1}^mc_{il}(x)c_{jl}(y)$.
 \end{itemize}
\medskip

Let $Y=(Y^1,\ldots,Y^N)$ be the solution to the stochastic
differential equation:
 \beqlb\label{2.0}
  \begin{cases}
    dY^1_t=b(Y^1_t)dt+e(Y^1_t)dB^1_t+c(Y^1_t)dW_t & \\
    \ldots\ldots & \\
    dY^N_t=b(Y^N_t)dt+e(Y^N_t)dB^N_t+c(Y^N_t)dW_t, &
  \end{cases}
 \eeqlb
where $W$ is an $\R^m$-valued Brownian motion, and $B^1,\ldots,B^N$
are mutually independent $E$-valued Brownian motion, which are
independent of $W$. Let $(S^N_t)_{t\geq 0}$ be the semigroup of the
diffusion $Y$ with generator $G_N$. Then for $f\in \DD(G_N)$, domain
of  $G_N$, it is easy to see that
   \beqnn
G_Nf(x_1,\ldots,x_N)&=&\sum_{p=1}^N\sum_{i=1}^db_i(x_p)\frac{\partial
f(x_1,\ldots,x_N)}{\partial
x_{p,i}}+\frac{1}{2}\sum_{p=1}^N\sum_{i,j=1}^dd_{ij}(x_p,x_p)
\frac{\partial^2f(x_1,\ldots,x_N)}{\partial x_{p,i}\partial x_{p,j}}\nonumber\\
& &+\frac{1}{2}\sum_{\gfc{p,q=1}{p\neq
q}}^N\sum_{i,j=1}^da_{ij}^{(m)}(x_p,x_q)\frac{\partial^2f(x_1,\ldots,x_N)}{\partial
x_{p,i}\partial x_{q,j}}.
 \eeqnn

We stress that under hypotheses (LU), the transition semigroup
$(S_t^N)_{t\geq 0}$ of the $Nd$-dimensional diffusion $Y$ has a
transition density, say $(p_t^N(x;y))_{t>0}$, and the semigroup is
both Feller and strong Feller; see, for instance, \cite[p.164]{RW00}
and \cite[p.227]{D65b}. Moreover, if hypotheses (LU) hold, then one
can modify the construction of \cite{SA01} to construct a dense
subset $D(\bar{E}^N)$ of $C_b(\bar{E}^N)$, satisfying
$D(\bar{E}^N)|_{E^N}\subset C^2_b(E^N)$, and extend $(S^N_t)_{t\geq
0}$ to a strongly continuous contraction semigroup
$(\bar{S}^N_t)_{t\geq 0}$ on $C_b(\bar{E}^N)$ such that
$D(\bar{E}^N)$ is invariant under $(\bar{S}^N_t)_{t\geq 0}$; see
\cite{D10}. Write $D(E^N)=D(\bar{E}^N)|_{E^N}$, class of functions
restricted to $E^N$. Note that functions in $D(E^N)$ are subject to
$\lim_{|x|\to\infty}G_Nf(x)=0$. Then $D(E):=D(E^1)\supset C^2_0(E)$
(space of functions together with their two continuous derivatives
vanishing at infinity).

To analyze $X^n=\{X_t^n:t\geq 0\}$, we associate to each $\alpha\in
I$ and $n$ a stopping time
 \beqnn \tau^{\alpha,n}=\left\{
  \begin{array}{ll}
    0 & \hbox{\quad if $\alpha_0>K_n$} \\
    \text{min}\{\frac{i+1}{n}: 0\leq i\leq |\alpha|, N^{\alpha|_i,n}=0\}
 & \hbox{\quad if this set is non-empty and $\alpha_0\leq K_n$} \\
    \frac{1+|\alpha|}{n} & \hbox{\quad otherwise}
  \end{array}
\right.
 \eeqnn
and define
 \beqnn
X^{\alpha,n}_t=\left\{
                \begin{array}{ll}
                  Y^{\alpha,n}_t & \hbox{\quad if $t<\tau^{\alpha,n}$} \\
                  \Delta & \hbox{\quad if $t\geq\tau^{\alpha,n}$.}
                \end{array}
              \right.
 \eeqnn
As a result, we have $X^n_t(h)=\frac{1}{n}\sum_{\alpha\sim_n
t}\hat{h}(X^{\alpha,n}_t)$ for measurable $h$, where $\hat{h}:=h$ on
$E$ and $\hat{h}(\Delta):=0$.

By It\^{o}'s formula it is easily verified that for each
$\alpha\sim_n k_n$, $t\in[k_n,k_n+a_n]$ and $f\in D(E)$
 \beqnn
M^{\alpha,k_n}_t(f)=
1_E(X^{\alpha,n}_{k_n})\left[f(Y^{\alpha,n}_t)-
f(Y^{\alpha,n}_{k_n})-\int^t_{k_n}Gf(Y^{\alpha,n}_u)du\right]
  \eeqnn
is an $(\FF^n_t)$-martingale with $G:=G_1$, moreover, from the
construction of $X^n$ we have
 \beqlb\label{2.1}
X^n_t(f)
&=&X^n_0(f)+M^{(n)}_t(f)+J^{(n)}_t(f)+N^{(n)}_t(f)+Z^{(n)}_t(f)+C^{(n)}_t(f)+H^{(n)}_t(f),
 \eeqlb
where
 \beqnn
M^{(n)}_t(f)&=&n^{-1}\sum_{r<k}\sum_{\alpha\sim_n
r_n}M^{\alpha,r_n}_{r_n+a_n}(f)[N^{\alpha,n}-\beta_n(Y^{\alpha,n}_{r_n+a_n})],\\
N^{(n)}_t(f)&=&n^{-1}\sum_{r<k}\sum_{\alpha\sim_n
r_n}\int^{r_n+a_n}_{r_n}\widehat{Gf}(X^{\alpha,n}_u)du[N^{\alpha,n}-\beta_n(Y^{\alpha,n}_{r_n+a_n})],
 \eeqnn
 \beqnn
 J^{(n)}_t(f)&=&n^{-1}\sum_{\alpha\sim_n
k_n}M^{\alpha,k_n}_t(f)+n^{-1}\sum_{r<k}\sum_{\alpha\sim_n
r_n}\int^{r_n+a_n}_{r_n}\widehat{Gf}(X^{\alpha,n}_u)du[\beta_n(Y^{\alpha,n}_{r_n+a_n})-1]\\
& &+n^{-1}\sum_{r<k}\sum_{\alpha\sim_n
r_n}\hat{f}(X^{\alpha,n}_{r_n})[\beta_n(Y^{\alpha,n}_{r_n+a_n})-\beta_n(Y^{\alpha,n}_{r_n})]\\
& &+n^{-1}\sum_{r<k}\sum_{\alpha\sim_n
r_n}M^{\alpha,r_n}_{r_n+a_n}(f)[\beta_n(Y^{\alpha,n}_{r_n+a_n})-\beta_n(Y^{\alpha,n}_{r_n})],
 \eeqnn
 \beqnn
Z^{(n)}_t(f)&=&n^{-1}\sum_{r<k}\sum_{\alpha\sim_n
r_n}\hat{f}(X^{\alpha,n}_{r_n})[N^{\alpha,n}-\beta_n(Y^{\alpha,n}_{r_n+a_n})]+n^{-1}\sum_{r<k}\sum_{\alpha\sim_n
r_n}M^{\alpha,r_n}_{r_n+a_n}(f)\beta_n(Y^{\alpha,n}_{r_n}),
 \eeqnn
 \beqnn
C^{(n)}_t(f)&=&n^{-1}\sum_{r<k}\sum_{\alpha\sim_n
r_n}\int^{r_n+a_n}_{r_n}
\widehat{Gf}(X^{\alpha,n}_u)du+n^{-1}\sum_{\alpha\sim_n k_n}
\int^t_{k_n}\widehat{Gf}(X^{\alpha,n}_u)du=\int^t_0X^n_u(Gf)du,
 \eeqnn
 \beqnn
H^{(n)}_t(f)&=&n^{-1}\sum_{r<k}\sum_{\alpha\sim_n
r_n}\hat{f}(X^{\alpha,n}_{r_n})[\beta_n(Y^{\alpha,n}_{r_n})-1]=\int^{k_n}_0X^n_{[ns]_n}(f\gamma_n)ds,
 \eeqnn
where $\hat{h}$ is defined as before. The major difference between
(\ref{2.1}) above and (A.3) of \cite{SA01} is the term $J^{(n)}$,
which is here more general.

It is well-known that to show the weak convergence of $\{X^n\}$
involves proving its tightness, deriving a martingale problem for
its limits and showing the uniqueness of solutions to the martingale
problem. Undoubtedly, the techniques of \cite{SA01} in deriving both
tightness and martingale characterizations are applicable here. For
uniqueness, the dual conditions of \cite{DK82} were used in
\cite{SA01}, while for the purpose of constructing a flow
superprocess with general branching variance and obtaining moment
formulas as well, the method of constructing directly dual processes
as in \cite{DLW01} is proved to be much more powerful.

 The lemma below shall play a fundamental role in proving
the tightness of $\{X^n\}$. Since the corresponding proof was not
given in \cite{SA01}, we provide one. Note that $X_0^n=
 \frac{1}{n}\sum^{K_n}_{i=1}\delta_{x^n_i}$.

\proclaim{\sc Lemma\n{2.1}.} {\sl Let $p$ be as in $(\ref{1.3})$ and
$T\geq 0$. If $X^n_0\Rightarrow\nu\in M_F(E)$, then
 \beqnn
C_T=\sup_{n\geq 1}\E\bigg(\sup_{0\leq t\leq
T}X^n_t(1)^2\bigg)<\infty \mbox{\ and \ }C'_T=\sup_{n\geq
1}\E\bigg(\sup_{0\leq t\leq T}X^n_t(1)^p\bigg)<\infty.
 \eeqnn}\endproclaim

\noindent{\it Proof.} Note that for non-negative Borel measurable
function $\phi$ on $E$ and for $t\in[k_n,k_n+a_n]$
 \beqlb\label{2.2}
\E X^n_t(\phi)\leq {\rm
e}^{[nt]_nC_\gamma}X_0^n(\E_\cdot[\phi(Y_t)]),
 \eeqlb
where $Y$ is as in (\ref{2.0}) with $N=1$, and $\E_y$ denotes the
conditional expectation given $Y_0=y$ (\cite[Lemma
II.3.3(a)]{P02}). It suffices to show that
  \beqlb\label{2.3}
C^{'}_T=\sup_{n\geq 1}\E\bigg(\sup_{0\leq t\leq
T}X^n_t(1)^p\bigg)<\infty.
 \eeqlb
In the following we let $C(u_1,\ldots,u_k)$ denote a constant
depending only on $u_1,\ldots,u_k$. Clearly
 \beqnn
\sup_{0\leq t\leq T}X^n_t(1)^p\leq
C(T,\gamma,p)\left(X^n_0(1)^p+\sup_{0\leq t\leq
T}|Z^{(n)}_t(1)|^p+\int^{[nT]_n}_0X^n_{[ns]_n}(1)^pds\right).
 \eeqnn
From $(\ref{1.2})$ and $(\ref{2.2})$, it follows that
$\{(Z^{(n)}_{k_n}(1),\FF^n_{k_n}): k=0,1,\ldots\}$ is a martingale.
Its predictable quadratic variation process  is calculated to be
 \beqnn
\<Z^{(n)}(1)\>_{k_n}:=\sum_{i=1}^k\E\left[\left(Z^{(n)}_{i_n}(1)-Z^{(n)}_{(i-1)_n}(1)\right)^2\bigg
|\FF^n_{(i-1)_n}\right]=\int^{k_n}_0X^n_{[ns]_n}(S^1_{a_n}(\sigma_n^2))ds.
 \eeqnn
Immediately
 \beqnn
\<Z^{(n)}(1)\>_{k_n}^{p/2}\leq
C(\sigma,p)\bigg(1+\int^{k_n}_0X^n_{[ns]_n}(1)^pds\bigg).
 \eeqnn
Then by Burkholder's inequality (\cite[p.152]{P02}) and H\"{o}lder's
inequality, we get
 \beqlb\label{2.4}
& &\E\bigg(\sup_{0\leq t\leq T}|Z^{(n)}_t(1)|^p\bigg)\\
&\leq&c\E\left[\left(\<Z^{(n)}(1)\>_{[nT]_n}\right)^{p/2}\right]\nonumber\\
& &+c\E\bigg[\bigg(\max_{0\leq k<[nT]}\bigg|n^{-1}\sum_{\alpha\sim_n
k_n}1_E(X^{\alpha,n}_{k_n})\left[N^{\alpha,n}
-\beta_n(Y^{\alpha,n}_{k_n+a_n})\right]\bigg|\bigg)^p\bigg]\nonumber\\
&\leq&\frac{c}{2}C_\sigma^pT^{p/(2q)}\bigg(1+\int^{[nT]_n}_0\E\left[X^n_{[ns]_n}(1)^p\right]ds\bigg)\nonumber\\
&
&+c\sum_{k=0}^{[nT]-1}\E\bigg[\E\bigg(\bigg|\frac{1}{n}\sum_{\alpha\sim_n
k_n}1_E(X^{\alpha,n}_{k_n})\left[N^{\alpha,n}-
\beta_n(Y^{\alpha,n}_{k_n+a_n})\right]\bigg|^p\bigg|\bar{\FF}^n_{k_n}\bigg)\bigg],\nonumber
 \eeqlb
where all the expectations above are allowed to be infinite, $c$ is
some constant depending only on $p$ and $1/p+1/q=1$. We shall use a
technique in \cite[pp.531-532]{SA01}. Fix $k$ and $n$. Note that
 \beqnn
\E\bigg(\bigg|n^{-1}\sum_{\alpha\sim_n
k_n}1_E(X^{\alpha,n}_{k_n})\left[N^{\alpha,n}-
\beta_n(Y^{\alpha,n}_{k_n+a_n})\right]\bigg|^p\bigg|\bar{\FF}^n_{k_n}\bigg)
=\E\bigg[\bigg|n^{-1}\sum_{i=1}^K\left[N^{\alpha^i,n}-
\beta_n(Y^{\alpha^i,n}_{k_n+a_n})\right]\bigg|^p\bigg],
 \eeqnn
where $K=nX^n_{k_n}(1)$, and $\alpha^1,\ldots,\alpha^K$ are the
labels of $K$ particles alive at time $k_n$. For $i=1,2,\ldots$, let
$\alpha^i\sim_n k_n$ such that $\alpha^i\neq\alpha^j$ if $i\neq j$
and define
 \beqnn
M_m(1):=n^{-1}\sum_{i=1}^m\left[N^{\alpha^i,n}-\beta_n(Y^{\alpha^i,n}_{k_n+a_n})\right]
\mbox{\ and\ }
\mathscr{G}_m:=\sigma(N^{\alpha^i,n}:i=1,\ldots,m)\bigvee\bar{\FF}^n_{k_n}.
 \eeqnn
Then clearly $Y^{\alpha^i,n}_{k_n}\in\mathscr{G}_m$ for all $i\geq
1$ and $\{(M_m(1),\mathscr{G}_m):m=1,2,\ldots\}$ is a square
integrable martingale by the fact that $N^{\alpha^m,n}$ is
conditionally independent of $\{N^{\alpha^i,n}:i=1,\ldots,m-1\}$
given $\bar{\FF}^n_{k_n}$. Similarly, we have
 \beqnn
\E\left(|M_m(1)|^p\right)&\leq&c\E\left(\<M(1)\>_m^{p/2}\right)+c\E\left(\max_{1\leq
l\leq m}|M_l(1)-M_{l-1}(1)|^p\right)\nonumber\\
&\leq&cn^{-p}C_\sigma^pm^{p/2}+c\sum_{l=1}^mn^{-p}\E\left(\left|N^{\alpha^l,n}
-\beta_n(Y^{\alpha^l,n}_{k_n+a_n})\right|^p\right)\leq
C(\gamma,\sigma,p,c)\left(\frac{\sqrt{m}}{n}\right)^p,
 \eeqnn
where $c$ is the same as the one in $(\ref{2.4})$. Recall that $p$
can be very close to 2. Therefore by the Cauchy-Schwarz inequality
and $C_T<\infty$ we obtain that
 \beqlb\label{2.5}
&&\sum_{k=0}^{[nT]-1}\E\bigg[\E\bigg(\bigg|\frac{1}{n}\sum_{\alpha\sim_n
k_n}1_E(X^{\alpha,n}_{k_n})\left[N^{\alpha,n}-
\beta_n(Y^{\alpha,n}_{k_n+a_n})\right]\bigg|^p\bigg|\bar{\FF}^n_{k_n}\bigg)\bigg]\\
&\leq&cC(\gamma,\sigma,p,c)\frac{[nT]}{n^{p/2}}\E\bigg[\bigg(\sup_{0\leq
t\leq T}X^n_t(1)\bigg)^{p/2}\bigg]\leq
cTC(\gamma,\sigma,p,c)n^{1-p/2}C_T^{p/4}.\nonumber
 \eeqlb
Now combine $(\ref{2.4})$ and $(\ref{2.5})$ to see that
 \beqnn
\E\bigg(\sup_{0\leq t\leq T}|Z^{(n)}_t(1)|^p\bigg)
&\leq&C(T,\gamma,\sigma,p,c)+C(T,\gamma,\sigma,p,c)\int^{[nT]_n}_0\E\left[X^n_{[ns]_n}(1)^p\right]ds,
 \eeqnn
 then $(\ref{2.3})$ follows by an analogous argument as in \cite[Lemma
II.4.6]{P02}. The proof is complete.

The following theorem is obviously an analogue of  a combination of
Propositions\n{A.3.10} and\n{A.3.12} and Lemma\n{A.3.13} of
\cite{SA01}, the proof of which certainly applies here except that
it suffices to prove (instead of the stronger square result)
 \beqnn
 \lim_{n\to\infty}\E\bigg(\sup_{0\leq t\leq
T}|J^{(n)}_t(f)|\bigg)=0
 \eeqnn
to obtain that $J^{(n)}(f)$ converges weakly to the zero process in
$D_\R(0,\infty)$. Now we state the result, and the interested reader
is referred to \cite{D101} for the detailed proof.

\proclaim{\sc Theorem\n{2.1}.} {\sl Suppose that the hypothesis (L)
holds. If $X^n_0\Rightarrow\nu\in M_F(E)$, then $\{X^n\}$ is tight
in $D_{M_F(E)}[0,\infty)$, and each limit point $X\in
C_{M_F(E)}[0,\infty)$ and is a solution to the following martingale
problem: For any $f\in D(E)$,
 \beqlb\label{2.6}
  Z_t(f)=X_t(f)-\nu(f)-\int_0^tX_s((G+\gamma)f)ds
 \eeqlb
 is a continuous square integrable martingale with $Z_0(f)=0$ and quadratic
 variation process
 \beqlb
 \label{2.61}
 \<Z(f)\>_t=\int_0^tX_s(\sigma^2f^2)ds+\int_0^t(X_s\times X_s)(\Lambda
 f)ds,\eeqlb
where $\Lambda
f(x,y)=\sum_{i,j=1}^da_{ij}^{(m)}(x,y)f'_i(x)f'_j(y)$.}\endproclaim

Then we shall prove the uniqueness of the martingale problem
(\ref{2.6}). Define for $F$ in some subset $\DD(\LL)$ (to be
specified) of the domain of an operator $\LL$ as follows
 \beqnn
 \LL F(\mu)\equiv(\LL
 F)(\mu)&:=&\int_{E}(G+\gamma)\left(\frac{dF(\mu)}{d\mu(x)}\right)\mu(dx)+\frac{1}{2}\int_{E}\sigma(x)^2\frac{d^2F(\mu)}
 {d\mu(x)^2}\mu(dx)\\
 & &+\frac{1}{2}\sum_{i,j=1}^d\int_{E}\int_{E}a_{ij}^{(m)}(x,y)
 \frac{\partial^2}{\partial x_i\partial
 y_j}\left(\frac{d^2F(\mu)}{d\mu(x)d\mu(y)}\right)\mu(dx)\mu(dy),
 \eeqnn
where $\frac{dF(\mu)}{d\mu(x)}:=\lim_{r\to
0+}\frac{1}{r}[F(\mu+r\delta_x)-F(\mu)], x\in E$, and similarly
$d^2F(\mu)/d\mu(x)d\mu(y)$ is
 defined with $F$ replaced by $d
 F(\mu)/d\mu(y)$. Let $X$ be a limit as in Theorem\n{2.1}. We will show that $X$
satisfies the martingale problem for $\LL$ and then construct the
dual process of $X$ to prove uniqueness.

Let $\DD(\LL)=\DD_1(\LL)\cup\DD_2(\LL)$, where $\DD_1(\LL)$ consists
of functions $F_f(\mu)=\<f,\mu^N\>$ with $f\in D(E^N)$, and
$\DD_2(\LL)$ denotes the class of functions
$F_{f,\phi}(\mu)=f(\mu(\phi_1),\ldots,\mu(\phi_N))$ with $f\in
C^2_b(\R^N)$ and $\phi=\{\phi_1,\ldots,\phi_N\}\subset D(E)$, and of
functions $F_{f,\phi}(\mu)=f(\mu(\phi))$ with $f\in
C^2_b([0,\infty))$ and $\phi\in D(E)^+$.  Let $\E_\nu$ denote the
conditional expectation given $X_0=\nu$.

\proclaim{\sc Lemma\n{2.2}.} {\sl $\E_\nu[X_t(1)^n]$ is locally
bounded in $t$ for each $n\geq 1$. Furthermore, $X$ is also a
solution to the martingale problem for $(\LL,\DD(\LL),\nu)$. That
is, for all $F\in\DD(\LL)$
 \beqlb\label{2.7}
F(X_t)-F(\nu)-\int_0^t\LL F(X_s)ds
 \eeqlb
is a continuous martingale with $X_0=\nu$.}
\endproclaim

\noindent{\it Proof.} Let $T_k=\inf\{t\geq 0: X_t(1)\geq k\}$. Then
$\{T_k\}$ is a non-decreasing sequence of stopping times. It is
easily seen from Lemma\n{2.1} that $T_k\to\infty$ as $k\to\infty$.
Fix $n\geq 1$. For each $k$, by It\^{o}'s formula we have
 \beqnn
X_{t\wedge T_k}(1)^n &=&X_0(1)^n+n\int^{t\wedge
T_k}_0X_s(1)^{n-1}dX_s(1)
     +\frac{n(n-1)}{2}\int^{t\wedge T_k}_0X_s(1)^{n-2}d\<Z(1)\>_s
 \eeqnn
 \beqnn
&=&X_0(1)^n+n\int^{t\wedge T_k}_0X_s(1)^{n-1}X_s(\gamma)ds
      +n\int^{t\wedge T_k}_0X_s(1)^{n-1}dZ_s(1)\\
& &+\frac{n(n-1)}{2}\int^{t\wedge T_k}_0X_s(1)^{n-2}X_s(\sigma^2)ds \\
&\leq&X_0(1)^n+nC_\gamma\int^{t\wedge
       T_k}_0X_s(1)^nds+\frac{n(n-1)}{2}C_\sigma^2\int^{t\wedge T_k}_0
            [1+X_s(1)^n]ds+\mbox{mart}.\nonumber
 \eeqnn
It follows that
 \beqnn
 \E[X_{t\wedge T_k}(1)^n]\leq
 X_0(1)^n+\frac{n(n-1)t}{2}C_\sigma^2+\left(nC_\gamma+
 \frac{n(n-1)}{2}C_\sigma^2\right)\int^t_0\E[X_{s\wedge
 T_k}(1)^n]ds.
 \eeqnn
An application of Gronwall's inequality and Fatou's lemma implies
the local boundedness of $\E[X_t(1)^n]$ in $t$. The martingale
property for $(\LL,\DD_1(\LL),\nu)$ is actually implied in
\cite[pp.537-539]{SA01}. It is sufficient to consider
$F_{f,\phi}=f(\mu(\phi_1),\ldots,\mu(\phi_N))$. Note that
 \beqnn
\LL
F_{f,\phi}(\mu)&=&\sum_{p=1}^Nf'_p(\mu(\phi_1),\ldots,\mu(\phi_N))\mu((G+\gamma)\phi_p)
+\frac{1}{2}\sum_{p,q=1}^Nf''_{pq}(\mu(\phi_1),\ldots,\mu(\phi_N))\mu(\sigma^2\phi_p\phi_q)\\
&
&+\frac{1}{2}\sum_{p,q=1}^Nf''_{pq}(\mu(\phi_1),\ldots,\mu(\phi_N))\int_E\int_E\sum_{i,j=1}^d
a^{(m)}_{ij}(x,y)\frac{\partial\phi_p(x)}{\partial
x_i}\frac{\partial\phi_q(y)}{\partial y_j}\mu(dx)\mu(dy).
 \eeqnn
By the martingale property for $(\LL,\DD_1(\LL),\nu)$ and It\^{o}'s
formula we have
 \beqnn
f(X_t(\phi_1),\ldots,X_t(\phi_N))&=&f(X_0(\phi_1),\ldots,X_0(\phi_N))+\mbox{mart.}\\
&
&+\int_0^t\sum_{p=1}^Nf'_p(X_s(\phi_1),\ldots,X_s(\phi_N))X_s((G+\gamma)\phi_p)ds\\
&
&+\int_0^t\frac{1}{2}\sum_{p,q=1}^Nf''_{pq}(X_s(\phi_1),\ldots,X_s(\phi_N))d\<Z(\phi_p),Z(\phi_q)\>_s.
 \eeqnn
Then the martingale property for $(\LL,\DD_2(\LL),\nu)$ follows once
we use polarization to see that
 \beqnn
\<Z(\phi_p),Z(\phi_q)\>_t=\int_0^tX_s(\sigma^2\phi_p\phi_q)ds+
\int_0^t\int_E\int_E\sum_{i,j=1}^da^{(m)}_{ij}(x,y)\frac{\partial\phi_p(x)}{\partial
x_i}\frac{\partial\phi_q(y)}{\partial y_j}X_s(dx)X_s(dy)ds.
 \eeqnn
This completes the proof.

It is worthwhile to notice that every solution to the martingale
problem (\ref{2.7}) is also such that (\ref{2.6}) is a continuous
local martingale with quadratic variation process given by
(\ref{2.61}); see, for instance, Theorem 4.8 of \cite{M96} and
Theorem 7.13 of \cite{L11}. Roughly speaking, the martingale
problems (\ref{2.6}) and (\ref{2.7}) are equivalent.

Before turning to our construction, observe that for $f\in D(E^N)$
 \beqlb\label{2.8}
\LL F_f(\mu)&=&F_{G_Nf}(\mu)+1/2\sum_{\gfc{p,q=1}{p\neq
q}}^N[F_{\Phi_{p,q}f}(\mu)-F_f(\mu)]\\
&
&+1/2\sum_{p=1}^N[F_{\Phi_pf}(\mu)-F_f(\mu)]+1/2N^2F_f(\mu)\nonumber\\
  &=&F_\mu(G_Nf,N)+1/2\sum_{\gfc{p,q=1}{p\neq
q}}^N[F_\mu(\Phi_{p,q}f,N-1)-F_\mu(f,N)]\nonumber \\
&
&+1/2\sum_{p=1}^N[F_\mu(\Phi_pf,N)-F_\mu(f,N)]+1/2N^2F_\mu(f,N),\nonumber
 \eeqlb
where for $h\in B(E^n)$ and $x=(x_1,\ldots,x_n)\in E^n$,
$F_\mu(h,n):= F_h(\mu)$, $\Phi_{p,q}: B(E^n)\to B(E^{n-1})$ is given
by
 \beqlb\label{2.9}
\Phi_{p,q}h(x_1,\ldots,x_{n-1}):=\sigma(x_{n-1})^2
h(x_1,\ldots,x_{n-1},\ldots,x_{n-1},\ldots,x_{n-2})
 \eeqlb
with $x_{n-1}$ in the positions of the $p$th and the $q$th variables
of $h$, and $\Phi_p: B(E^n)\to B(E^n)$ by
 \beqlb\label{2.10}
\Phi_ph(x):=2\gamma(x_p)h(x).
 \eeqlb

Based on (\ref{2.8}), we now construct a function-valued process of
$X$. Let $\mathbb{N}:=\{1,2,\ldots\}$. Let
$\mathbf{B}:=\cup_{n=0}^\infty B(E^n)$ be endowed with bounded
pointwise convergence on each $B(E^n)$, where $B(E^0):=\R$ and the
union is required to be disjoint union and so we do not view
$B(E^k)$ as a subset of $B(E^l)$ if $k<l$. Assume
$\{e_1,e_2,\ldots\}$ is a sequence of mutually independent unit
exponential random variables with $e_0:=0$. Define a sequence
$\Gamma=\{\Gamma_k:k=1,2,\ldots\}$ of random operators on
$\mathbf{B}$ and a $\mathbf{B}$-valued c\`{a}dl\`{a}g process
$L=\{L_t:t\geq 0\}$ as follows: Given a $\mathbf{B}$-valued random
variable $L_0$, independent of $\{e_1,e_2,\ldots\}$, define
recursively
 \beqnn
\left\{
       \begin{array}{ll}
      L_t=S^{N(L_{\tau_k})}_{t-\tau_k}\Gamma_kS^{N(L_{\tau_{k-1}})}_{\eta_k}\cdots
\Gamma_2S^{N(L_{\tau_1})}_{\eta_2}\Gamma_1S^{N(L_{\tau_0})}_{\eta_1}L_{\eta_0},
\quad\hbox{if\ } \tau_k\leq t<\tau_{k+1} \\
 \P\{\Gamma_{k+1}=\Phi_{p,q}|N(L_{\tau_k})=n_{k+1}\}=\P\{\Gamma_{k+1}=\Phi_p|N(L_{\tau_k})=n_{k+1}\}=
\frac{1}{n^2_{k+1}} \mbox{\ for\ } 1\leq p\neq q\leq n_{k+1}\\
L_{\tau_{k+1}}=\Gamma_{k+1}S^{N(L_{\tau_k})}_{\eta_{k+1}}\Gamma_kS^{N(L_{\tau_{k-1}})}_{\eta_k}\cdots
\Gamma_2S^{N(L_{\tau_1})}_{\eta_2}\Gamma_1S^{N(L_{\tau_0})}_{\eta_1}L_{\eta_0}, \quad k=0,1,2\ldots,\\
       \end{array}
     \right.\\
 \eeqnn
where $\eta_0=0, \eta_n=\frac{2e_n}{N(L_{\tau_{n-1}})^2}$,
$\tau_k=\sum_{i=0}^k\eta_i$ and $N(h):=l$ if $h\in B(E^l)$. Note
that given $L_0\in\mathbf{B}$, $\tau_k\to\infty$ almost surely as
$k\to\infty$ and thus $L_t$ is defined for all $t>0$. Set
$M_t=N(L_t)$. Then $(L,M)$ is a $\mathbf{B}\times\mathbb{N}$-valued
strong Markov process and shall serve as the dual process of $X$.
Let $\E_{h,n}$ denote the conditional expectation given
$(L_0,M_0)=(h,n)\in \mathbf{B}\times\mathbb{N}$ with $N(h)=n$ and
let $\mathscr{L}^\ast$ be the generator of $(L,M)$. Then from the
previous construction, one can verify, with elementary arguments,
that
 \beqlb\label{2.11}
&
&\E_{h,n}\left[\<L_t,\mu^{M_t}\>\exp\left\{\frac{1}{2}\int_0^tM_s^2ds\right\}\right]\\
&=&\<S_t^nh,\mu^n\>+\frac{1}{2}\sum_{\gfc{p,q=1}{p\neq
q}}^n\int_0^t\E_{\Phi_{p,q}S_s^nh,n-1}\left[\<L_{t-s},\mu^{M_{t-s}}\>
\exp\left\{\frac{1}{2}\int_0^{t-s}M_u^2du\right\}\right]ds\nonumber\\
&
&+\frac{1}{2}\sum_{p=1}^n\int_0^t\E_{\Phi_pS_s^nh,n}\left[\<L_{t-s},\mu^{M_{t-s}}\>
\exp\left\{\frac{1}{2}\int_0^{t-s}M_u^2du\right\}\right]ds\nonumber
 \eeqlb
and that $\LL^\ast F_{\mu}(f,N)=\LL F_f(\mu)-1/2N^2F_\mu(f,N)$ for
$f\in D(E^N)$.

\proclaim{\sc Theorem\n{2.2}.} {\sl Suppose that hypotheses $(LU)$
hold. Then for all $n\geq 1, t\geq 0$ and $h\in B(E^n)$ we have
 \beqlb\label{2.12}
\E\left[\<h,
X_t^n\>\right]=\E_{h,n}\left[\<L_t,\mu^{M_t}\>\exp\left\{\frac{1}{2}\int_0^tM_s^2ds\right\}\right],
 \eeqlb
where $X^n_t=X_t\times\cdots\times X_t\in M_F(E^n)$. Moreover,
uniqueness holds for the martingale problem $(\ref{2.7})$ and hence
for the martingale problem $(\ref{2.6})$.}
\endproclaim

\noindent{\it Proof.} In terms of Theorem\n{2.1}, Lemma\n{2.2} and
the relation (\ref{2.11}), the assertions follow in much the same
way as the proofs of Theorems 2.1 and 2.2 of \cite[p.7]{DLW01}. It
was pointed out in Remark 2.2 of He \cite{H09} that there is a gap
in the proof of Theorem 2.1 of \cite{DLW01} if $\sigma$ is just
bounded and measurable. This is because the function-valued process
$L$ is not always valued in the domain of $G_N$. However, if
$\sigma,\gamma,L_0\in C_b(E)$, then since the transition semigroup
of the underlying motion has regular transition density, hence
$L_t~(t\notin \{\tau_k\})$ does belong to the domain of $G_N$ and so
we can still use the associated martingale relation between $\LL$
and $\LL^*$. Consequently, the discussions in proving Theorem 2.1 of
\cite{DLW01} are applicable here since $\xi,\delta\in C_l(E)$. Note
that one may use Lemma\n{2.2} to verify that if $\gamma$ and
$\sigma$ are constants, then the total mass process $X(1)$ is a
diffusion process with generator $\gamma
x\frac{d}{dx}+\frac{\sigma^2x}{2}\frac{d^2}{dx^2}$ and $
 \E{\rm e}^{-\rho X_t(1)}={\rm e}^{-x\Psi_t(\rho)}$,
where $x=X_0(1)$ and $\Psi_t(\rho)=\frac{\rho{\rm e}^{\gamma
t}}{1+\frac{\sigma^2\rho({\rm e}^{\gamma t}-1)}{2\gamma}}$ (see
\cite[p.540]{SA01}).

It is natural to call an adapted c\`{a}dl\`{a}g process in $M_F(E)$
which satisfies the martingale problem (\ref{2.6}) a {\it
superprocess over a stochastic flow}, or simply {\it flow
superprocess $(G,\gamma,\sigma)$}.

In the remainder of this section, we shall derive the moment
formulas for $X$, the flow superprocess $(G,\gamma,\sigma)$ given by
Theorems\n{2.1} and\n{2.2}. Let $Y=(Y^1,\ldots,Y^N)$ be the $Nd$
dimensional diffusion process given by (\ref{2.0}), its semigroup
$S^N$. For $h\in B(E^N)$, define an operator $U^{(N)}$ by
 \beqnn
 U^{(N)}h=\frac{1}{2}\sum_{p\neq
 q\in\{1,\ldots,N\}}\Phi_{p,q}h,
  \eeqnn
 and a semigroup $T^N$ as follows
 \beqnn
 T_t^Nh(y)=\E_y\left[\exp\left\{\int_0^t\sum_{p=1}^N\gamma(Y^p(s))ds\right\}h(Y(t))\right].
 \eeqnn

 \proclaim{\sc Theorem\n{2.3}.} {\sl For $h\in B(E^n)$ and each
 $n\geq 1$
 \beqlb\label{2.13}
 \E_\nu\<h,X^n_t\>=\<T^n_th,\nu^n\>+\sum_{i=1}^{n-1}\<\int^t_0dt_1\int^{t_1}_0dt_2\cdots
 \int^{t_{i-1}}_0T^{n-i}_{t_i}\Pi^{(n)}(i-1;t)hdt_i,\nu^{n-i}\>,
 \eeqlb
where
$\Pi^{(n)}(i-1;t)=(U^{(n-(i-1))}T^{n-(i-1)}_{t_{i-1}-t_i}\cdots
 U^{(n-1)}T^{n-1}_{t_1-t_2})U^{(n)}T^n_{t-t_1}$ with
 $\Pi^{(n)}(0;t):=U^{(n)}T^n_{t-t_1}$, and $\int^t_0dt_1\int^{t_1}_0dt_2\cdots
 \int^{t_{i-1}}_0\cdot dt_{i-1}:=\int^t_0\cdot dt_1$ if $i=1$.
\endproclaim}

\noindent{\it Proof.} Let $h\in B(E^N)$. Define for $0\leq
t_{i+1}\leq t_i\leq\cdots\leq t_1\leq t$ the operators $V^{(N)}$ and
$\pi^{(N)}(i;t)$ respectively by
 \beqnn
 V^{(N)}h=\frac{1}{2}\sum_{p=1}^N\Phi_ph, \quad\mbox{ and }
  \pi^{(N)}(i;t)h=\left(V^{(N)}S^N_{t_i-t_{i+1}}\cdots
 V^{(N)}S^N_{t_1-t_2}\right)V^{(N)}S^N_{t-t_1}h
 \eeqnn
with $\pi^{(N)}(0;t):=V^{(N)}S^N_{t-t_1}$. Similarly
$\pi^{(N)}(k;s)$ is defined for $0\leq s_{k+1}\leq s_k\leq\cdots\leq
s_1\leq s$. To simplify notation, write
$V^{(N)}(x)=\frac{1}{2}\sum_{p=1}^N2\gamma(x_p)$. Then
$V^{(N)}h(x)=V^{(N)}(x)h(x)$. We first show that for any bounded
linear functionals $(V_t)_{t\geq 0}$ on $B(E^N)$
  \beqlb\label{2.14}
 &&\int^t_0dt_1\int^{t_1}_0dt_2\cdots\int^{t_i}_0V_{t_{i+1}}\left(S^N_{t_i-t_{i+1}}\pi^{(N)}(i-1;t)h\right)dt_{i+1}
 \\
 &=&\int^t_0V_{t_{i+1}}\bigg(\E_\cdot\bigg[\frac{1}{i!}\bigg(\int^t_{t_{i+1}}V^{(N)}(Y(s-{t_{i+1}}))ds\bigg)^i
 h(Y(t-{t_{i+1}}))\bigg]\bigg)dt_{i+1},\nonumber
 \eeqlb
 and
 \beqlb\label{2.141}
\int^t_0dt_1\int^{t_1}_0dt_2\cdots\int^{t_{i-1}}_0S^N_{t_i}\pi^{(N)}(i-1;t)hdt_i
=\E_\cdot\bigg[\frac{1}{i!}\bigg(\int^t_0V^{(N)}(Y(u))du\bigg)^i
 h(Y(t))\bigg].
 \eeqlb
We only consider (\ref{2.14}). Since
$S^N_{t_1-t_2}V^{(N)}S^N_{t-t_1}h(x)=\E_x[V^{(N)}(Y(t_1-t_2))h(Y(t-t_2))]$
by the Markov property of $Y$ as in (\ref{2.0}), so by Fubini's
theorem (\ref{2.14}) is clearly true for $i=1$. Suppose that it
holds for some $i\geq 1$. Then
 \beqnn
&&\int^t_0dt_1\int^{t_1}_0dt_2\cdots\int^{t_{i+1}}_0V_{t_{i+2}}\left(S^N_{t_{i+1}-t_{i+2}}\pi^{(N)}(i;t)h\right)dt_{i+2}\\
&=&\int^t_0dt_1\int^{t_1}_0V_{t_{i+2}}\bigg(\E_\cdot\bigg[\frac{1}{i!}
\bigg(\int^{t_1}_{t_{i+2}}V^{(N)}(Y(s-{t_{i+2}}))ds\bigg)^i
 (V^{(N)}S^N_{t-t_1}h)(Y(t_1-{t_{i+2}}))\bigg]\bigg)dt_{i+2}\\
&=&\int^t_0dt_{i+2}V_{t_{i+2}}\bigg(\E_\cdot
\bigg[\frac{1}{i!}\int^t_{t_{i+2}}\bigg(\int^{t_1}_{t_{i+2}}V^{(N)}(Y(s-{t_{i+2}}))ds\bigg)^i
 V^{(N)}(Y(t_1-t_{i+2}))dt_1h(Y(t-{t_{i+2}}))\bigg]\bigg)\\
&=&\int^t_0V_{t_{i+2}}\bigg(\E_\cdot
\bigg[\frac{1}{(i+1)!}\bigg(\int^t_{t_{i+2}}V^{(N)}(Y(s-{t_{i+2}}))ds\bigg)^{i+1}
h(Y(t-{t_{i+2}}))\bigg]\bigg)dt_{i+2}.
 \eeqnn
Thus (\ref{2.14}) is true for $i+1$ and hence for all $i\geq 1$.

By (\ref{2.11}), it is simple to see that (\ref{2.13}) holds for
$n=1$. Suppose that it is valid for some $n\geq 1$.  For $t\geq 0$,
write $H_t=\<L_t,\nu^{M_t}\>\exp\{\frac{1}{2}\int^t_0M_s^2ds\}$ and
$W^N_t(h)=\E_{h,N}[H_t]$ with $h\in B(E^N)$. Then $W^N_t$ is a
bounded linear functional on $B(E^N)$. In the remainder of this
proof, take $N=n+1$ and it is enough to consider $h\in
B(E^{n+1})^+$. By (\ref{2.11}), we have
 \beqnn
\E_{h,n+1}[H_s]&=&\<S^{n+1}_sh,\nu^{n+1}\>+\int^s_0
\E_{U^{(n+1)}S^{n+1}_{s-s_1}h,n}[H_{s_1}]ds_1+\int^s_0\E_{V^{(n+1)}S^{n+1}_{s-s_1}h,n+1}[H_{s_1}]ds_1\\
&=:&W_s(h)+\int^s_0\E_{V^{(n+1)}S^{n+1}_{s-s_1}h,n+1}[H_{s_1}]ds_1,\quad
s\geq 0.
 \eeqnn
Make repeated use of the above relation to conclude that
 \beqlb\label{2.15}
\E_{h,n+1}[H_s]&=&W_s(h)+\sum_{i=1}^k\int^s_0ds_1\int^{s_1}_0ds_2\cdots\int^{s_{i-1}}_0
W_{s_i}\left(\pi^{(n+1)}(i-1;s)h\right)ds_i\\
 &
 &+\int^s_0ds_1\int^{s_1}_0ds_2\cdots\int^{s_k}_0\E_{\pi^{(n+1)}(k;s)h,n+1}[H_{s_{k+1}}]ds_{k+1},\quad
 k\geq 1.\nonumber
 \eeqlb
Since by (\ref{2.12}) as well as the induction assumption
 \beqnn
W_{s_i}\left(\pi^{(n+1)}(i-1;s)h\right)&=&\<S^{n+1}_{s_i}\left(\pi^{(n+1)}(i-1;s)h\right),\nu^{n+1}\>\\
 & &+\int^{s_i}_0W^n_u\left(U^{(n+1)}S^{(n+1)}_{s_i-u}(\pi^{(n+1)}(i-1;s)h)\right)du,
 \eeqnn
hence in terms of (\ref{2.14}) and (\ref{2.141}), we see that
 \beqnn
&&\sum_{i=1}^k\int^s_0ds_1\int^{s_1}_0ds_2\cdots\int^{s_{i-1}}_0W_{s_i}
\left(\pi^{(n+1)}(i-1;s)h\right)ds_i\\
&=&\<\sum_{i=1}^k\int^s_0ds_1\int^{s_1}_0ds_2\cdots\int^{s_{i-1}}_0
S^{n+1}_{s_i}\left(\pi^{(n+1)}(i-1;s)h\right)ds_i,\nu^{n+1}\>\\
& &+\sum_{i=1}^k\int^s_0ds_1\int^{s_1}_0ds_2\cdots\int^{s_i}_0
W^n_{s_{i+1}}\left(U^{(n+1)}S^{(n+1)}_{s_i-s_{i+1}}(\pi^{(n+1)}(i-1;s)h)\right)ds_{i+1}\\
&=&\<\sum_{i=1}^k\E_\cdot\bigg[\frac{1}{i!}\left(\int^s_0V^{(n+1)}(Y(u))du\right)^i
h(Y(s))\bigg],\nu^{n+1}\>\\
& &+\sum_{i=1}^k\int^s_0W^n_{s_{i+1}}\bigg(U^{(n+1)}\bigg(\E_\cdot
\bigg[\frac{1}{i!}\bigg(\int^s_{s_{i+1}}V^{(n+1)}(Y_{s_i-s_{i+1}})ds_i\bigg)^i
h(Y(s-s_{i+1}))\bigg]\bigg)\bigg)ds_{i+1}
 \eeqnn
 \beqnn
&\underrightarrow{k\to\infty}&
\<\E_\cdot\bigg[\sum_{i=1}^\infty\frac{1}{i!}\left(\int^s_0V^{(n+1)}(Y(u))du\right)^i
h(Y(s))\bigg],\nu^{n+1}\>\\
 & &+\int^s_0W^n_{s_{i+1}}\bigg(U^{(n+1)}\bigg(\E_\cdot
\bigg[\sum_{i=1}^\infty\frac{1}{i!}\bigg(\int^s_{s_{i+1}}V^{(n+1)}(Y(s_i-s_{i+1}))ds_i\bigg)^i
h(Y(s-s_{i+1}))\bigg]\bigg)\bigg)ds_{i+1}.
 \eeqnn
Note that the last term in (\ref{2.15}) tends to zero as
$k\to\infty$ for $X$ has locally bounded moments of any order. Then
letting $k\to\infty$ we get
 \beqnn
\E_{h,n+1}[H_s]&=&\<S^{n+1}_sh,\nu^{n+1}\>
+\<\E_\cdot\bigg[\sum_{i=1}^\infty\frac{1}{i!}\left(\int^s_0V^{(n+1)}(Y(u))du\right)^i
h(Y(s))\bigg],\nu^{n+1}\>\\
& &+\int^s_0W^n_{s_1}\left(U^{(n+1)}S^{n+1}_{s-s_1}h\right)ds_1\\
& &+\int^s_0W^n_{s_1}\bigg(U^{(n+1)}\bigg(\E_\cdot
\bigg[\sum_{i=1}^\infty\frac{1}{i!}\bigg(\int^s_{s_1}V^{(n+1)}(Y_{u-s_1})du\bigg)^i
h(Y(s-s_1))\bigg]\bigg)\bigg)ds_1\\
&=&\<T^{n+1}_sh,\nu^{n+1}\>+\int^s_0W^n_{s_1}\left(U^{(n+1)}T^{n+1}_{s-s_1}h\right)ds_1.
 \eeqnn
Now a simple variable change together with (\ref{2.12}) implies that
(\ref{2.13}) holds for $n+1$. This completes the proof.

The formulas below were established in \cite{SA01} under the
condition of binary branching, using branching particle systems
approximation. They are immediate from Theorem\n{2.3} and the Markov
property of $X$.

 \proclaim{\sc
Corollary\n{2.1}.} {\sl If $\gamma$ and $\sigma$ are constants, then
for $h,h_1,h_2\in B(E)$ and $0\leq s\leq t$
 \beqnn
\E_\nu\left[X_t(h)\right]&=&{\rm e}^{\gamma t}\nu(S_t^1h),
 \eeqnn
 and
 \beqnn \E_\nu\left[
X_s(h_1)X_t(h_2)\right]&=&{\rm
e}^{\gamma(s+t)}\<S_s^2(h_1S_{t-s}^1h_2),\nu^2\>+\sigma^2{\rm
e}^{\gamma(s+t)}\int_0^s{\rm e}^{-\gamma
u}\<S_u^1[S_{s-u}^2(h_1S_{t-s}^1h_2)],\nu\>du,
 \eeqnn
 where $
S_u^1[S_{s-u}^2(h_1S_{t-s}^1h_2)](x)
:=\int_E\int_E\int_Eh_1(w)S^1_{t-s}h_2(z)S_{s-u}^2(y,y;dw,dz)S_u^1(x,dy)$.
}\endproclaim

We note that by constructing  a stochastic integral similarly as in
\cite{DLW01}, it is quite simple to obtain the first moment of $X$
while it does not seem obvious to derive higher moments.

\medskip

{\bf3. Measurable spatially dependent branching.} In this part we
shall construct, via approximation, flow superprocesses with drift
function $\gamma\in B(E)$ and branching variance $\sigma^2\in
B(E)^+$. To this aim, choose functions $\{\gamma^{(i)}\}\subset
C_l(E)$ and
 $\{\sigma^{(i)}\}\subset C_l(E)^+$ such that
$\gamma^{(i)}(x)\to\gamma(x), \sigma^{(i)}(x)\to\sigma(x)$ as
$i\to\infty$ for $\lambda$-a.e. $x\in E$. Here $\lambda$ is Lebesgue
measure on $E$. Let $\{X^{(i)}\}$ be the flow superprocesses
$(G,\gamma^{(i)},\sigma^{(i)})$ given by Theorems\n{2.1} and\n{2.2}.
Proceed as in the previous section to construct a function-valued
process $L=\{L_t:t\geq 0\}$ based on $(\gamma,\sigma)$. Define
naturally mappings $\Phi^{(i)}_{p,q}$ and $\Phi^{(i)}_p$ as in
(\ref{2.9}) and (\ref{2.10}) with $\gamma$ and $\sigma$ replaced by
$\gamma^{(i)}$ and $\sigma^{(i)}$, respectively. In an obvious way
we can construct operators $\Gamma^{(i)}=\{\Gamma^{(i)}_k\}$ and
function-valued processes $L^{(i)}=\{L_t^{(i)}:t\geq 0\}$. Define
$M_t^{(i)}=N(L_t^{(i)})$, which is clearly independent of $i$ and
hence write $M_t^{(i)}=:M_t$.

 \proclaim{\sc Lemma\n{3.1}.} {\sl Suppose that $\mu_i\Rightarrow \mu$ in $M_F(E)$.
 Then for $t\geq 0$ and $h\in \mathbf{B}$ with $N(h)=n$
  \beqnn
  \lim_{i\to\infty}\E_{h,n}\left[\<L_t^{(i)},\mu_i^{M_t}\>\exp\left\{\frac{1}{2}\int_0^tM_s^2ds\right\}\right]
  =\E_{h,n}\left[\<L_t,\mu^{M_t}\>\exp\left\{\frac{1}{2}\int_0^tM_s^2ds\right\}\right].
  \eeqnn }\endproclaim

\noindent{\it Proof.} Let $L_{(h,n)}=\{L_{(h,n)}(t):t\geq0\}$ denote
the process $L$ with $L_0=h$ and $N(h)=n$. Let
$L^{(i)}_{(h,n)}=\{L^{(i)}_{(h,n)}(t):t\geq 0\}$ stand for the
process $L^{(i)}$ with initial value $L_0^{(i)}=h$ and $N(h)=n$. By
(\ref{2.11}), we obtain that
 \beqlb\label{3.1}
&&\E_{h,n}\left[\<L_t^{(i)},\mu_i^{M_t}\>\exp\left\{\frac{1}{2}\int_0^tM_s^2ds\right\}\right]\\
&=&\E\left[\<L^{(i)}_{(h,n)}(t),\mu_i^{M_t}\>\exp\left\{\frac{1}{2}\int_0^tM_s^2ds\right\}\right]\nonumber\\
&=&\<S_t^nh,\mu_i^n\>+\E\left[\int_0^t\<L^{(i)}_{(\Phi^{(2,n)}_{t-s}h,n-1)}(s),\mu_i^{M_s}\>
\exp\left\{\frac{1}{2}\int_0^sM_u^2du\right\}ds\right]\nonumber\\
&
&+\E\left[\int_0^t\<L^{(i)}_{(\Phi^{(1,n)}_{t-s}h,n)}(s),\mu_i^{M_s}\>
\exp\left\{\frac{1}{2}\int_0^sM_u^2du\right\}ds\right],\nonumber
 \eeqlb
where $\Phi^{(2,n)}_s=\frac{1}{2}\sum_{p\neq
q\in\{1,\ldots,n\}}\Phi_{p,q}^{(i)}S_s^n$, and
$\Phi^{(1,n)}_s=\frac{1}{2}\sum_{p=1}^n\Phi_p^{(i)}S_s^n$. For
notational simplicity, write
$n_k=N(L^{(i)}_{\tau_k})=N(L_{\tau_k})$. Since $(S^N_t)_{t\geq 0}$
has transition density $(p^N_t(x;y))_{t>0}$ and since
$\gamma^{(i)}(x)\to\gamma(x)$ and $\sigma^{(i)}(x)\to\sigma(x)$
$\lambda$-a.e. $x$, we can use dominated convergence and induction
to see that for $\tau_k<t<\tau_{k+1}$ with $k=0,1,\ldots$
 \beqnn
 \left(S^{n_k}_{t-\tau_k}\Gamma_k^{(i)}S^{n_{k-1}}_{\eta_k}\cdots
\Gamma^{(i)}_1S^{n_0}_{\eta_1}L^{(i)}_{\eta_0}\right)(z)
\longrightarrow
\left(S^{n_k}_{t-\tau_k}\Gamma_kS^{n_{k-1}}_{\eta_k}\cdots
\Gamma_1S^{n_0}_{\eta_1}L_{\eta_0}\right)(z)
 \eeqnn
for all $z\in E^{n_k}$,  and
 \beqnn
 \left(\Gamma^{(i)}_{k+1}S^{n_k}_{t-\tau_k}\Gamma_k^{(i)}S^{n_{k-1}}_{\eta_k}\cdots
\Gamma^{(i)}_1S^{n_0}_{\eta_1}L^{(i)}_{\eta_0}\right)(z)
\longrightarrow
\left(\Gamma_{k+1}S^{n_k}_{t-\tau_k}\Gamma_kS^{n_{k-1}}_{\eta_k}\cdots
\Gamma_1S^{n_0}_{\eta_1}L_{\eta_0}\right)(z)
 \eeqnn
 for $\lambda_{k+1}$-a.e. $z\in E^{n_{k+1}}$, where $\lambda_{k+1}$
denotes Lebesgue measure on $E^{n_{k+1}}$. Now fix an arbitrary $k$.
Then for $g\in C_b(E^{n_k}\times E^{n_k})$ with compact support, by
Fubini's theorem and dominated convergence again we have
 \beqnn
&&\int_{E^{n_k}}\int_{E^{n_k}}g(x,y)\left(\Gamma_k^{(i)}S^{n_{k-1}}_{\eta_k}\cdots
\Gamma^{(i)}_1S^{n_0}_{\eta_1}L^{(i)}_{\eta_0}\right)(y)p^{n_k}(x;y)
dy\mu_i^{n_k}(dx)\\
&=&\int_{E^{n_k}}\int_{E^{n_k}}g(x,y)p^{n_k}(x;y)\mu_i^{n_k}(dx)
\left(\Gamma_k^{(i)}S^{n_{k-1}}_{\eta_k}\cdots
\Gamma^{(i)}_1S^{n_0}_{\eta_1}L^{(i)}_{\eta_0}\right)(y)dy\\
&\to&\int_{E^{n_k}}\int_{E^{n_k}}g(x,y)p^{n_k}(x;y)\mu^{n_k}(dx)
\left(\Gamma_kS^{n_{k-1}}_{\eta_k}\cdots
\Gamma_1S^{n_0}_{\eta_1}L_{\eta_0}\right)(y)dy.
 \eeqnn
By \cite[Proposition 4.4, p.112]{EK86}, for each $k$ the measures
$\{(\Gamma_k^{(i)}S^{n_{k-1}}_{\eta_k}\cdots
\Gamma^{(i)}_1S^{n_0}_{\eta_1}L^{(i)}_{\eta_0})(y)p^{n_k}(x;y) $
$dy\mu_i^{n_k}(dx)\}$ on $E^{n_k}\times E^{n_k}$ are weakly
convergent. It is then evident that the required result follows from
(\ref{2.11}) and (\ref{3.1}) and hence the proof is finished.

The following theorem asserts the existence of flow superprocesses
with bounded branching variance and drift function.

\proclaim{\sc Theorem\n{3.1}. }  {\sl $X^{(i)}\Rightarrow X $ in
$C_{M_F(E)}[0,\infty)$ and $X$ satisfies the martingale problem
(\ref{2.7}).}
 \endproclaim

\noindent{\it Proof.} The proof can be proceeded as the one of
Theorem\n{5.2} of \cite{DLW01} and hence is omitted.

\medskip

\noindent{\it Remark 3.1.} It is clear that the moments of flow
superprocess $(G,\gamma,\sigma)$, $X$, are still given by
Theorem\n{2.3}. Moreover, both its mean measure $m_1(A):=\E X_t(A)$
and covariance measure $m_2(A\times B):=\E(X_t(A)X_t(B))$ are
absolutely continuous with respect to Lebesgue measure since
$\gamma$ is bounded and the semigroups $S^1$ and $S^2$ have
densities.

\medskip
 \textbf{Acknowledgments.} I would like to thank Professor
Zenghu Li for his supervision.
\medskip

\center{REFERENCES}

\end{document}